\documentclass{article}
\usepackage[english]{babel}
\usepackage[utf8]{inputenc}
\usepackage{fancyhdr}
\usepackage[numbers]{natbib}

\usepackage{graphicx}
\usepackage{amsmath}
\usepackage{relsize}
\usepackage{multicol}
\usepackage{ragged2e}
\usepackage{fontenc}
\usepackage{mathcomp}
\usepackage{latexsym}
\usepackage{amssymb}
\usepackage{times}
\newtheorem{Teorema}{Theorem}[section]

\newtheorem{Definicija}[Teorema]{Definition}
\newtheorem{Posledica}[Teorema]{Corollary}

\newtheorem{Lema}[Teorema]{Lemma}
\newtheorem{Primedba}[Teorema]{Remark}
\numberwithin{equation}{section}
\usepackage{dsfont}
\newcommand\scalemath[2]{\scalebox{#1}{\mbox{\ensuremath{\displaystyle #2}}}} % paket za smanjivanje matrice
\numberwithin{equation}{section}
\hyphenchar\font=-1
\title{Perturbing the spectrum of operator $T_n^d(A)$}

\author{Nikola Sarajlija\footnote{nikola.sarajlija@dmi.uns.ac.rs, University of Novi Sad, Faculty of Sciences, Novi Sad, 21000, Serbia}} 
\begin{document}
\maketitle

\begin{abstract}
Let $T_n^d(A)$ denote a partial upper triangular operator matrix whose diagonal entries are given and the others unknown. In this article we have aim to find characterizations of (left,right) invertibility of $T_n^d(A)$ in terms of diagonal entries solely, and hence we provide statements which generalize and correct results of Zhang S., Wu Z. (2012). We pose our results without invoking separability condition, thus improving results of Zhang S., Wu Z. (2012), and we give appropriate n-dimensional analogues, without assuming separability as well. We recover many perturbation results of Djordjevi\'c D. S. (2002), and obtain some results of Du H. K., Pan J. (1994) and Han J. K., Lee H. Y., Lee W. Y. (2000) in the case of the Hilbert space setting.
\end{abstract}



\section{Introduction and preliminaries}

Spectral properties of upper triangular block operators are thoroughly studied by numerous authors (see \cite{CHEN}, \cite{OPERATORTHEORY}, \cite{KOLUNDZIJA2}, \cite{KOLUNDZIJA3}, \cite{HWANG}, \cite{KOLUNDZIJA}, \cite{KINEZI}, \cite{SARAJLIJA2}, \cite{WU}). The former are often encountered while investigating spectral properties of operators acting on a sum of Banach or Hilbert spaces. However, specialists have usually examined the case of $2\times2$ block operators, while the case of general $n\times n$ block operators has been neglected. This article deals with the study of this latter case. Furthermore, the usual framework in studying $2\times2$ block operators has been the separable Hilbert space setting, while our results hold without assuming separability.%Namely, let $T$ acts on a direct topological sum of Banach spaces $X$ and $Y$, $T:X\oplus Y\rightarrow X\oplus Y$, and assume that $X$ is invariant for $T$, $T(X)\subseteq X$. Then $T$ admits an upper triangular block representation $T=\begin{bmatrix}T_1 & T_2\\ 0 & T_3\end{bmatrix}:\begin{bmatrix}X\\ Y\end{bmatrix}\rightarrow\begin{bmatrix}X\\ Y\end{bmatrix}$. We can also consider $T$ acting on a direct topological sum of an arbitrary number of Banach spaces, which yields upper triangular block operators of a higher dimension . Not many articles have been published in a connection with upper triangular block operators of dimension $n>2$. This article is ought to be an exception.

We consider partial upper triangular operator matrices of arbitrary dimension $n\geq 2$. By "partial" we understand a matrix whose some of the entries are unknown, while the others are given. We investigate some completion problems related to (left, right) invertibility of such matrices, thus obtaining some per\-tur\-bation results as consequences.  Invertibility of $2\times 2$ block matrices of general form has been discussed in \cite{CHEN} by Chen and Hai, and M. Kolund\v{z}ija extended their result to case of Banach spaces \cite{KOLUNDZIJA}. Invertibility of $2\times 2$ block operators of upper triangular form has been discussed in \cite{KINEZI}, \cite{ZHANG}, and we extend results from \cite{ZHANG} to arbitrary $n\times n$ upper triangular operators. Similar investigations related to Fredholmness and Weylness of these block operators have been studied in \cite{SARAJLIJA2}.  Notice that there are various results related to different types of spectra (\cite{HUANG}, \cite{KINEZI2}, \cite{WU}, \cite{ZGUITTI}) of block operators, but there are only a few results treating the usual (left,right) spectrum.  

The article is organized as follows. In the rest of this section we give notations and some basics on invertibility theory with a few auxiliary results. Afterwards, Section 2 deals with characterizing (left, right) invertibility  of upper triangular operator matrices in the setting of arbitrary Banach spaces.

Let $ X, X_1,..., X_n$ be complex Banach spaces. Notation $\mathcal{B}( X_i, X_j)$ stands for the collection of all linear and bounded operators from $ X_i$ to $ X_j$, where we put $\mathcal{B}( X):=\mathcal{B}( X, X)$. If $T\in\mathcal{B}( X_i, X_j)$, $\mathcal{N}(T)$ and $\mathcal{R}(T)$ denote the kernel and the range of $T$, respectively. It is well known that $\mathcal{N}(T)$ is closed.

Let $T\in\mathcal{B}( X)$. Then $T$ is left (right) invertible if $\widetilde{T}T=I$ ($T\widetilde{T}=I$) for some $\widetilde{T}\in\mathcal{B}( X)$. $T$ is invertible if $T$ is both left and right invertible. Families of left (right) invertible operators on $ X$ are denoted by $\mathcal{G}_l( X)$ ($\mathcal{G}_r( X)$). Hence, $\mathcal{G}( X)=\mathcal{G}_l( X)\cap\mathcal{G}_r( X)$ is the family of invertible operators on $ X$.

Let $D_1\in\mathcal{B}( X_1),\ D_2\in\mathcal{B}( X_2),...,D_n\in\mathcal{B}( X_n)$ be given. Partial upper triangular operator matrix of dimension $n$ is denoted by
\begin{equation}\label{OSNOVNI}
T_n^d(A)=
\scalemath{0.9}{
\begin{bmatrix} 
    D_1 & A_{12} & A_{13} & ... & A_{1,n-1} & A_{1n}\\
    0 & D_2 & A_{23} & ... & A_{2,n-1} & A_{2n}\\
    0 &  0 & D_3 & ... & A_{3,n-1} & A_{3n}\\
    \vdots & \vdots & \vdots & \ddots & \vdots & \vdots\\
    0 & 0 & 0 & ... & D_{n-1} & A_{n-1,n}\\
    0 & 0 & 0 & ... & 0 & D_n      
\end{bmatrix}}\in\mathcal{B}( X_1\oplus  X_2\oplus\cdots\oplus  X_n),
\end{equation}
where $A:=(A_{12},\ A_{13},...,\ A_{ij},...,\ A_{n-1,n})$ is an operator tuple consisting of unknown variables $A_{ij}\in\mathcal{B}( X_j, X_i)$, $1\leq i<j\leq n,\ n\geq2$. Let $\mathcal{B}_n$ stands for the collection of all described tuples $A=(A_{ij})$.  This notation originated in \cite{WU}, and the present author used it in \cite{SARAJLIJA2}. We emphasize that entries of $T_n^d(A)$ are linear and bounded operators acting on appropriate Banach spaces. It is possible to replace the setting of Banach spaces with Hilbert $C^*$-modules \cite{STEFAN}, but we do not pursue this point any further.

In the sequel, we try to find an answer to the following question.

\noindent\textbf{Question 1.} Can we find a characterization for (left, right) invertibility of $T_n^d(A)$ in terms of (left, right) invertibility of its diagonal entries $D_i$?

If it is possible, than we immediately have an answer to

\noindent\textbf{Question 2.} What can we say about $\bigcap\limits_{A\in\mathcal{B}_n}\sigma_*(T_n^d(A))$, when $\sigma_*$ runs through set $\lbrace\sigma_{l},\sigma_{r},\sigma\rbrace$?

We shall present all of our results using some notions from Fredholm theory. Name\-ly, we will put into use notions of nullity and deficiency of operator. Precise definition follows \cite{ZANA}. Let operator $T\in\mathcal{B}(X)$. Then its nullity and deficiency are the following quantities, respectively: $\alpha(T)=\dim\mathcal{N}(T)$ and $\beta(T)=\dim X/\mathcal{R}(T)$. These quantities are nonnegative integers or $+\infty$. The following lemma enlights the reason for using such terminology in this article.
\begin{Lema}\label{VEZA3}
Let $T\in\mathcal{B}(X)$. The following equivalences hold:
$$\mathit{T\ is\ left\ invertible\ \Leftrightarrow\ \alpha(T)=0\ and\ \mathcal{R}(T)\ is\ closed\ and\ complemented;}$$
$$\mathit{T\ is\ right\ invertible\ \Leftrightarrow\ \beta(T)=0\ and\ \mathcal{N}(T)\ is\ complemented.}$$
\end{Lema}
If $X$ is a Hilbert space, then every closed subspace of $X$ is complemented. In that case Lemma \ref{VEZA3} becomes simpler.

Corresponding spectra of an operator $T\in\mathcal{B}( X)$ are defined as follows: the left spectrum is $\sigma_{l}(T)=\lbrace\lambda\in\mathds{C}: \lambda-T\not\in\mathcal{G}_{l}( X)\rbrace$, the right spectrum is $\sigma_{r}(T)=\lbrace\lambda\in\mathds{C}: \lambda-T\not\in\mathcal{G}_{r}( X)\rbrace$,
while the spectrum of $T$ is the following set: $\sigma(T)=\lbrace\lambda\in\mathds{C}: \lambda-T\not\in\mathcal{G}( X)\rbrace$. All of these spectra are compact nonempty subsets of the complex plane. We write $\rho_{l}(T), \rho_{r}(T), \rho(T)$ for their complements, respectively.

The following result is simple and well-known.

\begin{Lema}\label{POMOCNALEMA2}
Let $T_n^d(A)\in\mathcal{B}(X_1\oplus\cdots\oplus X_n).$ Then:\\
$(i)$ $\sigma_{l}(D_1)\subseteq\sigma_{l}(T_n^d(A))$;\\
$(ii)$ $\sigma_{r}(D_n)\subseteq\sigma_{r}(T_n^d(A)).$
\end{Lema}

\begin{Primedba}
If $Y$ is a complemented subspace of $X$, we denote a topological complement of $Y$ in $X$ by $Y_1$. Specially, in the rest of this work, if $Y=\mathcal{N}(T)$ ($Y=\overline{\mathcal{R}(T)}$), we use $\mathcal{M}(T)$ ($\mathcal{K}(T)$) to denote its topological complement in $X$. 
\end{Primedba}

One important class of operators whose properties we aim to exploit is the class of inner regular operators. $T\in\mathcal{B}(X)$ is inner regular if $T=T\widetilde{T}T$ for some $\widetilde{T}\in\mathcal{B}(X)$. For brevity, such $T$ will be called regular. By \cite[Corollary 1.1.5]{DJORDJEVIC}, $T$ is regular if and only if $\mathcal{N}(T)$, $\mathcal{R}(T)$ are closed and complemented subspaces. It is obvious that operators in each of the classes $\mathcal{G}_l(X),\mathcal{G}_r(X),\mathcal{G}(X)$ are (inner) regular. 

We will use the following notion that is due to D. S. Djordjevi\'{c} \cite{OPERATORTHEORY}.

\begin{Definicija}
We say that a Banach space $X$ can be embedded in a Banach space $Y$ if there is a left invertible operator from $X$ to $Y$. In that case we write $X\preceq Y$. Embedding $X\preceq Y$ is essential if $Y/\mathcal{R}(T)$ is an infinite dimensional linear space for every $T\in\mathcal{B}(X,Y)$.  In that case we write $X\prec Y$.
\end{Definicija}
\begin{Primedba}
It is obvious that $X\preceq Y$ if and only if there exists a right invertible operator from $Y$ to $X$.\\
If $X$ and $Y$ are Hilbert spaces, then $X\preceq Y$ $\Leftrightarrow$ $\dim X\leq \dim Y$, and $X\prec Y$ $\Leftrightarrow$ $\dim X<\dim Y$ and $\dim Y=\infty$. Here, $\dim$ stands for the orthogonal dimension.
\end{Primedba}

\section{Invertible completions of $T_n^d(A)$}
We start with a result which deals with left invertibility of $T_n^d(A)$. 

\begin{Teorema}\label{LEVIINVERZ'}
Let $D_1\in\mathcal{B}(X_1),\ D_2\in\mathcal{B}(X_2),...,D_n\in\mathcal{B}(X_n)$. Assume that $D_s,\ 2\leq s\leq n-1$, are regular operators. Consider the following statements:\\
$(i)$ $(a)$ $D_1\in\mathcal{G}_l(X_1)$;\\
\hspace*{5.5mm}$(b)$ $D_n$ is regular and $\mathcal{N}(D_i)\preceq \mathcal{K}(D_{i-1})$ for every $2\leq i\leq n$;\\[3mm]
$(ii)$ There exists $A\in\mathcal{B}_n$ such that $T_n^d(A)\in\mathcal{G}_l(X_1\oplus\cdots\oplus X_n)$;\\[3mm]
$(iii)$ $(a)$ $D_1\in\mathcal{G}_l(X_1)$;\\
\hspace*{8.2mm}$(b)$ $\bigoplus\limits_{s=1}^{i-1} \mathcal{K}(D_s)\prec\mathcal{N}(D_i)$ does not hold for $2\leq i\leq n$.\\[3mm]
Then $(i) \Rightarrow (ii).$\\
If $X_1,...,X_n$ are infinite dimensional Hilbert spaces, then $(ii)\Rightarrow(iii)$.
\end{Teorema}
\textbf{Proof. }$(i)\Rightarrow(ii)$

In this case it holds  $\alpha(D_1)=0$, $\mathcal{R}(D_s)$ is closed for all $1\leq s\leq n$ and  $\mathcal{N}(D_i)\preceq \mathcal{K}(D_{i-1})$ for every $2\leq i\leq n$. By Lemma \ref{VEZA3} we need to find $A\in\mathcal{B}_n$ such that $\alpha(T_n^d(A))=0$ and $\mathcal{R}(T_n^d(A))$ is closed and complemented.  We choose $A=(A_{ij})_{1\leq i<j\leq n}$ so that $A_{ij}=0$ if $j-i\neq 1$, that is we place all nonzero operators of tuple $A$ on the superdiagonal. It remains to define $A_{ij}$ for $j=i+1$, $1\leq i<n$. First notice that $A_{i,i+1}: X_{i+1}\rightarrow  X_i$. Since all of diagonal entries have closed range, we know that $ X_{i+1}=\mathcal{N}(D_{i+1})\oplus \mathcal{M}(D_{i+1})$, $ X_i=\mathcal{K}(D_i)\oplus\mathcal{R}(D_i)$, and we have  $\mathcal{N}(D_{i+1})\preceq \mathcal{K}(D_i)$. It follows that there is a left invertible operator $J_{i}:\mathcal{N}(D_{i+1})\rightarrow \mathcal{K}(D_i)$. We put $A_{i,i+1}=\begin{bmatrix}J_{i} & 0\\ 0 & 0\end{bmatrix}:\begin{bmatrix}\mathcal{N}(D_{i+1})\\ \mathcal{M}(D_{i+1})\end{bmatrix}\rightarrow\begin{bmatrix}\mathcal{K}(D_i)\\ \mathcal{R}(D_i)\end{bmatrix}$, and we implement this procedure for all $1\leq i\leq n-1$. Notice that $\mathcal{R}(A_{i,i+1})$ is contained in a subspace which is complementary to $\mathcal{R}(D_i)$ for each $1\leq i\leq n-1.$

Now we have chosen our $A$, we show that $\mathcal{N}(T_n^d(A))\cong\mathcal{N}(D_1)$, implying \\$\alpha(T_n^d(A))=\alpha(D_1)=0$. Let us put $T_n^d(A)x=0$, where $x=x_1+\cdots+x_n\in  X_1\oplus\cdots\oplus  X_n$. The previous equality is then equivalent to the following system of equations
$$
\begin{bmatrix}D_1x_1+A_{12}x_2\\ D_2x_2+A_{23}x_3\\ \vdots\\ D_{n-1}x_{n-1}+A_{n-1,n}x_n\\ D_nx_n\end{bmatrix}=\begin{bmatrix}0\\ 0\\ \vdots\\ 0\\ 0\end{bmatrix}.
$$ 
Last equation gives $x_n\in\mathcal{N}(D_n)$. Since $\mathcal{R}(A_{s,s+1})$ is contained in a subspace which is complementary to $\mathcal{R}(D_s)$ for all $1\leq s\leq n-1$, we have $D_sx_s=A_{s,s+1}x_{x+1}=0$ for all $1\leq s\leq n-1$. That is, $x_i\in\mathcal{N}(D_i)$ for every $1\leq i\leq n$, and $J_sx_{s+1}=0$ for every $1\leq s\leq n-1.$ Due to left invertibility of $J_s$ we get $x_s=0$ for $2\leq s\leq n$, which proves the claim. Therefore, $\alpha(T_n^d(A))=\alpha(D_1)=0$.

Next, we show that $\mathcal{R}(T_n^d(A))$ is closed and complemented. Left invertibility of $J_i$'s implies the existence of closed subspaces $U_i$ of $\mathcal{K}(D_i)$ such that $\mathcal{K}(D_i)=\mathcal{R}(J_i)\oplus U_i$, $1\leq i\leq n-1$ (Lemma \ref{VEZA3}). It means that $ X_1\oplus  X_2\oplus\cdots\oplus  X_n=\mathcal{R}(D_1)\oplus\mathcal{R}(J_{1})\oplus U_1\oplus\mathcal{R}(D_2)\oplus\mathcal{R}(J_{2})\oplus U_2\oplus\cdots\oplus\mathcal{R}(D_{n-1})\oplus\mathcal{R}(J_{n-1})\oplus U_{n-1}\oplus\mathcal{R}(D_n)\oplus\mathcal{K}(D_n)$. It is not hard to see that $\mathcal{R}(T_n^d(A))=$\\$\mathcal{R}(D_1)\oplus\mathcal{R}(J_{1})\oplus\mathcal{R}(D_2)\oplus\mathcal{R}(J_{2})\oplus\cdots\oplus\mathcal{R}(D_{n-1})\oplus\mathcal{R}(J_{n-1})\oplus\mathcal{R}(D_n)$. Comparing these equalities, one easily sees that $\mathcal{R}(T_n^d(A))$ is closed and complemented (see also \cite[Theorem 3.6]{RAKIC}).

\noindent$(ii)\Rightarrow (iii)$\\[3mm] 
\hspace*{6mm}Assume that $T_n^d(A)$ is left invertible and $X_1,...,X_n$ are Hilbert spaces. Then $D_1\in\mathcal{G}_l( X_1)$ (Lemma \ref{POMOCNALEMA2}). Assume that $(iii)(b)$ fails. Then there exists some $j\in\lbrace 2,..., n\rbrace$ such that $\alpha(D_j)>\sum\limits_{s=1}^{j-1}\beta(D_s)$.

We use a method similar to that in \cite{SARAJLIJA2}. We know that for each $A\in\mathcal{B}_n$, the operator matrix $T_n^d(A)$ as an operator from $ X_1\oplus\mathcal{N}(D_2)^\bot\oplus\mathcal{N}(D_2)\oplus\mathcal{N}(D_3)^\bot\oplus\mathcal{N}(D_3)\oplus\cdots\oplus\mathcal{N}(D_n)^\bot\oplus\mathcal{N}(D_n)$ into $\mathcal{R}(D_1)\oplus\mathcal{R}(D_1)^\bot\oplus\mathcal{R}(D_2)\oplus\mathcal{R}(D_2)^\bot\oplus\cdots\oplus\mathcal{R}(D_{n-1})\oplus\mathcal{R}(D_{n-1})^\bot\oplus X_n$ admits the following block representation
\begin{equation}\label{MATRICA}
T_n^d(A)=\scalemath{0.85}
{\begin{bmatrix} 
    D_1^{(1)} & A_{12}^{(1)} & A_{12}^{(2)} & A_{13}^{(1)} & A_{13}^{(2)} & ... & A_{1n}^{(1)} & A_{1n}^{(2)}\\
    0 & A_{12}^{(3)} & A_{12}^{(4)} & A_{13}^{(3)} & A_{13}^{(4)} & ... & A_{1n}^{(3)} & A_{1n}^{(4)}\\
    0 & D_2^{(1)} & 0 & A_{23}^{(1)} & A_{23}^{(2)} & ... & A_{2n}^{(1)} & A_{2n}^{(2)}\\
    0 & 0 & 0 & A_{23}^{(3)} & A_{23}^{(4)} & ... & A_{2n}^{(3)} & A_{2n}^{(4)}\\
    0 & 0 & 0 & D_{3}^{(1)} & 0 & ... & A_{3n}^{(1)} & A_{3n}^{(2)}\\
    0 & 0 & 0 & 0 & 0 & ... & A_{3n}^{(3)} & A_{3n}^{(4)}\\
    \vdots & \vdots & \vdots & \vdots & \vdots & \ddots & \vdots & \vdots\\
    0 & 0 & 0 & 0 & 0 & ... & A_{n-1,n}^{(1)} & A_{n-1,n}^{(2)}\\
    0 & 0 & 0 & 0 & 0 & ... & A_{n-1,n}^{(3)} & A_{n-1,n}^{(4)}\\
    0 & 0 & 0 & 0 & 0 & ... & D_n^{(1)} & 0\\
\end{bmatrix}}.
\end{equation}
Notice that $D_s^{(1)}$, $1\leq s\leq n-1$ are invertible, and $D_n^{(1)}$ is injective. Therefore, there exist invertible operator matrices $U$ and $V$ such that 
\begin{equation}\label{MATRICA2}
UT_n^d(A)V=\scalemath{0.85}{\begin{bmatrix}
D_1^{(1)} & 0 & 0 & 0 & 0 & ... & 0 & 0\\
0 & 0 & A_{12}^{(4)} & 0 & A_{13}^{(4)} & ... & A_{1n}^{(3)} & A_{1n}^{(4)}\\
0 & D_2^{(1)} & 0 & 0 & 0 & ... & 0 & 0\\
0 & 0 & 0 & 0 & A_{23}^{(4)} & ... & A_{2n}^{(3)} & A_{2n}^{(4)}\\
0 & 0 & 0 & D_3^{(1)} & 0 & ... & 0 & 0\\
0 & 0 & 0 & 0 & 0 & ... & A_{3n}^{(3)} & A_{3n}^{(4)}\\
\vdots & \vdots & \vdots & \vdots & \vdots & \ddots & \vdots & \vdots\\
0 & 0 & 0 & 0 & 0 & ... & 0 & 0\\
0 & 0 & 0 & 0 & 0 & ... & A_{n-1,n}^{(3)} & A_{n-1,n}^{(4)}\\
0 & 0 & 0 & 0 & 0 & ... & D_n^{(1)} & 0\\
\end{bmatrix}
}
\end{equation}
We will explain the construction of matrices $U$ and $V$ in more details. It is known that elementary transformations of a matrix can be carried out by multiplying the matrix with elementary matrices. In that way, since $D_1^{(1)}$, $D_2^{(1)}$,...$D_{n-1}^{(1)}$ are invertible, by multiplying the matrix $T_n^d(A)$ with suitable elementary matrices from the left, we ,,destroy" operators $A_{ij}^{(1)}$ and $A_{ij}^{(3)}$, where $1\leq i,j\leq n-1$. The product of those matrices is our matrix U. Now, analogously, we multiply $T_n^d(A)$ with suitable elementary matrices from the right in order to ,,destroy'' operators $A_{ij}^{(2)}$; the product of those matrices equals matrix $V$.\\
Note that $A_{ij}^{(3)}$ and $A_{ij}^{(4)}$ in (\ref{MATRICA2}) are not the original ones from (\ref{MATRICA}) in general, but we still use them for convenience. Now, it is obvious that if (\ref{MATRICA2}) is left invertible, then since $D_n^{(1)}$ is injective,
\begin{equation}\label{MATRICA3}
\begin{bmatrix}
A_{12}^{(4)} & A_{13}^{(4)} & A_{14}^{(4)} & ... & A_{1n}^{(4)}\\
0          & A_{23}^{(4)} & A_{24}^{(4)} & ... & A_{2n}^{(4)}\\
0           & 0         & A_{34}^{(4)} & ... & A_{3n}^{(4)}\\
\vdots   &   \vdots & \vdots & \ddots & \vdots\\
0          &   0        &    0      & ... & A_{n-1,n}^{(4)}
\end{bmatrix}
:
\begin{bmatrix}
\mathcal{N}(D_2)\\
\mathcal{N}(D_3)\\
\mathcal{N}(D_4)\\
\vdots\\
\mathcal{N}(D_n)
\end{bmatrix}
\rightarrow
\begin{bmatrix}
\mathcal{R}(D_1)^\bot\\
\mathcal{R}(D_2)^\bot\\
\mathcal{R}(D_3)^\bot\\
\vdots\\
\mathcal{R}(D_{n-1})^\bot
\end{bmatrix}
\end{equation}
is injective. Since $\alpha(D_j)>\sum\limits_{s=1}^{j-1}\beta(D_s)$ it follows that 
$$
\begin{bmatrix}
A_{1j}^{(4)}\\
A_{2j}^{(4)}\\
A_{3j}^{(4)}\\
\vdots\\
A_{j-1,j}^{(4)}
\end{bmatrix}:\mathcal{N}(D_j)\rightarrow\begin{bmatrix}\mathcal{R}(D_1)^\bot\\\mathcal{R}(D_2)^\bot\\ \vdots \\\mathcal{R}(D_{j-1})^\bot\end{bmatrix}
$$
is not injective, and hence operator defined in (\ref{MATRICA3}) is not injective for every $A\in\mathcal{B}_n$. Contradiction. This proves the desired. $\square$

\begin{Primedba}\label{VAZISVUDA}
Notice the validity of Theorem \ref{LEVIINVERZ'} without assuming separability of $ X_1,..., X_n$.
\end{Primedba}
\begin{Posledica}\label{POSLEDICA'}
Let $D_1\in\mathcal{B}(X_1),\ D_2\in\mathcal{B}(X_2),...,D_n\in\mathcal{B}(X_n)$. Assume that $D_s-\lambda,\ 2\leq s\leq n-1,\ \lambda\in\mathds{C}$ are regular operators. Then
$$
\bigcap\limits_{A\in\mathcal{B}_n}\sigma_{l}(T_n^d(A))\subseteq\sigma_{l}(D_1)\cup\Big(\bigcup\limits_{k=2}^{n}\Delta_k'\Big)\cup\Delta'',
$$
where
$$
\Delta_k':=\Big\lbrace\lambda\in\mathds{C}:\ \mathcal{N}(D_k-\lambda)\preceq \mathcal{K}(D_{k-1}-\lambda)\ does\ not\ hold\Big\rbrace,\ 2\leq k\leq n,
$$
$$
\Delta''=\Big\lbrace\lambda\in\mathds{C}:\ D_n-\lambda\ is\ not\ regular\Big\rbrace.
$$
If $X_1,...,X_n$ are infinite dimensional Hilbert spaces, then
$$
\sigma_{l}(D_1)\cup\Big(\bigcup\limits_{k=2}^{n}\Delta_k\Big)\subseteq\bigcap\limits_{A\in\mathcal{B}_n}\sigma_{l}(T_n^d(A)),
$$
where 
$$
\Delta_k:=\Big\lbrace\lambda\in\mathds{C}:\ \bigoplus\limits_{s=1}^{k-1}\mathcal{K}(D_s-\lambda)\prec\mathcal{N}(D_k-\lambda)\ holds\Big\rbrace,\ 2\leq k\leq n.
$$
\end{Posledica}
\begin{Primedba}
Obviously, $\Delta_k\subseteq\Delta_k'$ for $2\leq k\leq n$ holds.
\end{Primedba}
If $n=2$, we recover a result from \cite{OPERATORTHEORY}.
\begin{Teorema}(\cite[Theorem 5.2]{OPERATORTHEORY})\label{LEVIINVERZn=2'}
Let $D_1\in\mathcal{B}(X_1), D_2\in\mathcal{B}(X_2)$. Consider the following statements:\\[1mm]
$(i)$ $(a)$ $D_1\in\mathcal{G}_l(X_1)$;\\
\hspace*{5.5mm}$(b)$ $D_2$ is regular;\\
\hspace*{5.5mm}$(c)$ $\mathcal{N}(D_2)\preceq \mathcal{K}(D_1)$;\\[1mm]
$(ii)$ There exists $A\in\mathcal{B}_2$ such that $T_2^d(A)\in\mathcal{G}_l(X_1\oplus X_2)$;\\[1mm]
$(iii)$ $(a)$ $D_1\in\mathcal{G}_l(X_1)$;\\
\hspace*{8.2mm}$(b)$ $\mathcal{K}(D_1)\prec\mathcal{N}(D_2)$ does not hold.\\[3mm]
Then $(i)\Rightarrow(ii)$.\\
If $X_1,X_2$ are infinite dimensional Hilbert spaces, then $(ii)\Rightarrow(iii)$. 
\end{Teorema}
\begin{Posledica}(\cite[Corollary 5.3]{OPERATORTHEORY})\label{POSLEDICA2'}
Let $D_1\in\mathcal{B}(X_1),\ D_2\in\mathcal{B}(X_2)$. Then
$$
\bigcap\limits_{A\in\mathcal{B}_2}\sigma_{l}(T_2^d(A))\subseteq\sigma_{l}(D_1)\cup\Delta_2'\cup\Delta'',
$$
where
$$
\Delta_2':=\Big\lbrace\lambda\in\mathds{C}:\ \mathcal{N}(D_2-\lambda)\preceq \mathcal{K}(D_1-\lambda)\ does\ not\ hold\Big\rbrace,
$$
$$
\Delta''=\Big\lbrace\lambda\in\mathds{C}:\ D_2-\lambda\ is\ not\ regular\Big\rbrace.
$$
If $X_1,X_2$ are infinite dimensional Hilbert spaces, then 
$$
\sigma_{l}(D_1)\cup\Delta_2\subseteq\bigcap\limits_{A\in\mathcal{B}_2}\sigma_{l}(T_2^d(A)),
$$
where 
$$
\Delta_2:=\Big\lbrace\lambda\in\mathds{C}:\ \mathcal{K}(D_1-\lambda)\prec\mathcal{N}(D_2-\lambda)\ holds\Big\rbrace.
$$
\end{Posledica}
\begin{Primedba}
One might conjecture that the left invertible $T_2^d(A)$ must have $D_2$ with closed range. However, this is not the case. See \cite[Lemma 2]{WU} and \cite[Example 3]{HWANG}.
\end{Primedba}
Notice that Theorem \ref{LEVIINVERZn=2'} is a correct version of \cite[Theorem 2.1]{ZHANG}. There are several remarks concerning Theorem 2.1 in \cite{ZHANG}. First of all, in the notation of \cite{ZHANG}, condition $(i)(b)$ of Theorem \ref{LEVIINVERZn=2'} is omitted in \cite[Theorem 2.1]{ZHANG}, which is an oversight. Without that condition direction $(ii)\Rightarrow(i)$ in \cite[Theorem 2.1]{ZHANG} need not hold. Namely, the choice of $Q$ in the proof of part $(ii)\Rightarrow(iv)$ implies $\mathcal{R}(M_Q)= X\oplus\mathcal{R}(B)$, and for $\mathcal{R}(M_Q)$ to be closed (Lemma \ref{VEZA3}) we must assume that $\mathcal{R}(B)$ is closed. Furthermore, if $\mathcal{R}(B)$ is closed, notice that condition ($\beta(A)=\infty$ or ($B\in\Phi_+(\mathcal{K})$ and $\alpha(B)\leq\beta(A)$)) in \cite[Theorem 2.1]{ZHANG} is equivalent to a simple condition $\alpha(B)\leq\beta(A)$, which is condition $(i)(c)$ in Theorem \ref{LEVIINVERZn=2'} interpreted in the setting of Hilbert spaces. \\
Similar reasoning holds for \cite[Theorem 2.2]{ZHANG}.

Now, we provide results dealing with right invertibility of $T_n^d(A)$.
\begin{Teorema}\label{DESNIINVERZ'}
Let $D_1\in\mathcal{B}(X_1),\ D_2\in\mathcal{B}(X_2),...,D_n\in\mathcal{B}(X_n)$. Assume that $D_s,\ 2\leq s\leq n-1$ are regular operators. Consider the following statements:\\
$(i)$   $(a)$ $D_n\in\mathcal{G}_r(X_n)$;\\
\hspace*{5.5mm}$(b)$ $D_1$ is regular and $\mathcal{K}(D_i)\preceq\mathcal{N}(D_{i+1})$ for every $1\leq i\leq n-1$;\\[3mm]
$(ii)$ There exists $A\in\mathcal{B}_n$ such that $T_n^d(A)\in\mathcal{G}_r(X_1\oplus\cdots\oplus X_n)$;\\[3mm]
$(iii)$ $(a)$ $D_n\in\mathcal{G}_r(X_n)$;\\
\hspace*{8.2mm}$(b)$ $\bigoplus\limits_{s=i+1}^n\mathcal{N}(D_s)\prec \mathcal{K}(D_i)$ does not hold for $1\leq i\leq n-1$.\\[3mm]
Then $(i) \Rightarrow (ii)$.\\
If $X_1,...,X_n$ are infinite dimensional Hilbert spaces, then $(ii)\Rightarrow(iii)$.
\end{Teorema}
\textbf{Proof. }$(i)\Rightarrow(ii)$

In this case it holds  $\beta(D_n)=0$, $\mathcal{R}(D_s)$ is closed for all $1\leq s\leq n$ and $\mathcal{K}(D_i)\preceq\mathcal{N}(D_{i+1})$ for every $1\leq i\leq n-1$. By Lemma \ref{VEZA3}, we need to find $A\in\mathcal{B}_n$ such that $\beta(T_n^d(A))=0$ and $\mathcal{N}(T_n^d(A))$ is closed and complemented.  We choose $A=(A_{ij})_{1\leq i<j\leq n}$ so that $A_{ij}=0$ if $j-i\neq 1$, that is we place all nonzero operators of tuple $A$ on the superdiagonal. It remains to define $A_{ij}$ for $j=i+1$, $1\leq i< n$. First notice that $A_{i,i+1}: X_{i+1}\rightarrow  X_i$. Since all of diagonal entries have closed and complemented range and kernel, we know that $ X_{i+1}=\mathcal{N}(D_{i+1})\oplus \mathcal{M}(D_{i+1})$, $ X_i=\mathcal{K}(D_i)\oplus\mathcal{R}(D_i)$, and we have  $\mathcal{K}(D_i)\preceq\mathcal{N}(D_{i+1})$. It follows that there is a right invertible operator $J_{i}:\mathcal{N}(D_{i+1})\rightarrow \mathcal{K}(D_i)$. We put $A_{i,i+1}=\begin{bmatrix}J_{i} & 0\\ 0 & 0\end{bmatrix}:\begin{bmatrix}\mathcal{N}(D_{i+1})\\ \mathcal{M}(D_{i+1})\end{bmatrix}\rightarrow\begin{bmatrix}\mathcal{K}(D_i)\\ \mathcal{R}(D_i)\end{bmatrix}$, and we implement this procedure for all $1\leq i\leq n-1$. 

Notice that $\mathcal{R}(A_{i,i+1})=\mathcal{K}(D_i)$ for each $1\leq i\leq n-1.$ Therefore, it is immediate that $\mathcal{R}(T_n^d(A))=\mathcal{R}(D_1)\oplus\mathcal{R}(A_{12})\oplus\mathcal{R}(D_2)\oplus\mathcal{R}(A_{23})\oplus\cdots\oplus\mathcal{R}(D_{n-1})\oplus\mathcal{R}(A_{n-1,n})\oplus\mathcal{R}(D_n)$ is equal to $X_1\oplus\cdots\oplus X_n$, that is $T_n^d(A)$ is surjective.

Now we show that $T_n^d(A)$ has a complemented kernel. First, by Lemma \ref{VEZA3}, there exist closed subspaces $V_{i+1}$ of $\mathcal{N}(D_{i+1})$ such that $\mathcal{N}(D_{i+1})=\mathcal{N}(J_i)\oplus V_{i+1}$, $1\leq i\leq n-1$. It means that $ X_1\oplus  X_2\oplus\cdots\oplus  X_n=\mathcal{N}(D_1)\oplus\mathcal{M}(D_1)\oplus\mathcal{N}(D_2)\oplus\mathcal{N}(J_{1})\oplus V_2\oplus\cdots\oplus\mathcal{N}(D_{n})\oplus\mathcal{N}(J_{n-1})\oplus V_{n}$. Second, direct computation shows that $\mathcal{N}(T_n^d(A))\cong\mathcal{N}(D_1)\oplus\mathcal{N}(J_1)\oplus\cdots\oplus\mathcal{N}(J_{n-1})$. Comparing these equalities, and consulting Theorem 3.6 from \cite{RAKIC}, we conclude that $\mathcal{N}(T_n^d(A))$ is closed and complemented.

$(ii)\Rightarrow(iii)$

This implication follows directly from part $(ii)\Rightarrow(iii)$ of Theorem \ref{LEVIINVERZ'} by employing dual relations $\mathcal{N}(T)=\mathcal{R}(T^*)^\bot$, $\mathcal{N}(T^*)=\mathcal{R}(T)^\bot$. $\square$

\begin{Posledica}\label{POSLEDICA3'}
Let $D_1\in\mathcal{B}(X_1),\ D_2\in\mathcal{B}(X_2),...,D_n\in\mathcal{B}(X_n)$. Assume that $D_s-\lambda,\ 2\leq s\leq n-1,\ \lambda\in\mathds{C}$ are regular operators. Then
$$
\bigcap\limits_{A\in\mathcal{B}_n}\sigma_{r}(T_n^d(A))\subseteq\\\sigma_{r}(D_n)\cup\Big(\bigcup\limits_{k=2}^{n-1}\Delta_k'\Big)\cup\Delta'',
$$
where
$$
\Delta_k':=\Big\lbrace\lambda\in\mathds{C}:\ \mathcal{K}(D_k-\lambda)\preceq\mathcal{N}(D_{k+1}-\lambda)\ does\ not\ hold\Big\rbrace,\ 1\leq k\leq n-1,
$$
$$
\Delta'':=\Big\lbrace\lambda\in\mathds{C}:\ D_1-\lambda\ is\ not\ regular\Big\rbrace.
$$
If $X_1,...,X_n$ are infinite dimensional Hilbert spaces, then
$$
\sigma_{r}(D_n)\cup\Big(\bigcup\limits_{k=1}^{n-1}\Delta_k\Big)\subseteq\\\bigcap\limits_{A\in\mathcal{B}_n}\sigma_{r}(T_n^d(A)),
$$
where
$$
\Delta_k=\Big\lbrace\lambda\in\mathds{C}:\ \bigoplus\limits_{s=k+1}^n\mathcal{N}(D_s-\lambda)\prec \mathcal{K}(D_k-\lambda)\ holds\Big\rbrace,\ 1\leq k\leq n-1.
$$
\end{Posledica}
\begin{Primedba}
Obviously, $\Delta_k\subseteq\Delta_k'$ for $1\leq k\leq n-1$ holds.
\end{Primedba}

If $n=2$, we recover more results from \cite{OPERATORTHEORY}.

\begin{Teorema}(\cite[Theorem 5.4]{OPERATORTHEORY})\label{DESNIINVERZn=2'}
Let $D_1\in\mathcal{B}(X_1), D_2\in\mathcal{B}(X_2)$. Consider the following statements:\\[1mm]
$(i)$   $(a)$ $D_2\in\mathcal{G}_r(X_2)$;\\
\hspace*{5.5mm}$(b)$ $D_1$ is regular;\\
\hspace*{5.5mm}$(c)$ $\mathcal{K}(D_1)\preceq\mathcal{N}(D_{2})$;\\[3mm]
$(ii)$ There exists $A\in\mathcal{B}_2$ such that $T_2^d(A)\in\mathcal{G}_r(X_1\oplus X_2)$;\\[3mm]
$(iii)$ $(a)$ $D_2\in\mathcal{G}_r(X_2)$;\\
\hspace*{8.2mm}$(b)$ $\mathcal{N}(D_2)\prec \mathcal{K}(D_1)$ does not hold.\\[3mm]
Then $(i) \Rightarrow (ii)$.\\
If $X_1,X_2$ are infinite dimensional Hilbert spaces, then $(ii)\Rightarrow(iii)$.
\end{Teorema}
\begin{Posledica}(\cite[Corollary 5.5]{OPERATORTHEORY})\label{POSLEDICA4'}
Let $D_1\in\mathcal{B}(X_1),\ D_2\in\mathcal{B}(X_2)$. Then
$$
\bigcap\limits_{A\in\mathcal{B}_2}\sigma_{r}(T_2^d(A))\subseteq\\\sigma_{r}(D_2)\cup\Delta_1'\cup\Delta'',
$$
where
$$
\Delta_1':=\Big\lbrace\lambda\in\mathds{C}:\ \mathcal{K}(D_1-\lambda)\preceq\mathcal{N}(D_{2}-\lambda)\ does\ not\ hold\Big\rbrace,
$$
$$
\Delta'':=\Big\lbrace\lambda\in\mathds{C}:\ D_1-\lambda\ is\ not\ regular \Big\rbrace.
$$
If $X_1,X_2$ are infinite dimensional Hilbert spaces, then
$$
\sigma_{r}(D_2)\cup\Delta_1\subseteq\\\bigcap\limits_{A\in\mathcal{B}_2}\sigma_{r}(T_2^d(A)),
$$
where
$$
\Delta_1:=\Big\lbrace\lambda\in\mathds{C}:\ \mathcal{N}(D_2-\lambda)\prec \mathcal{K}(D_1-\lambda)\ holds\Big\rbrace.
$$

\end{Posledica}

We finish our investigations with results regarding invertibility of $T_n^d(A)$.

\begin{Teorema}\label{INVERZ'}
Let $D_1\in\mathcal{B}(X_1),\ D_2\in\mathcal{B}(X_2),...,D_n\in\mathcal{B}(X_n)$. Assume that all $D_s,\ 2\leq s\leq n-1$, are inner regular operators. Consider the following statements:\\
$(i)$   $(a)$ $D_1\in\mathcal{G}_l(X_1)$ and $D_n\in\mathcal{G}_r(X_n)$;\\
\hspace*{5.5mm}$(b)$ $\mathcal{N}(D_{i+1})\cong \mathcal{K}(D_i)$ for $1\leq i\leq n-1$;\\[3mm]
$(ii)$ There exists $A\in\mathcal{B}_n$ such that $T_n^d(A)\in\mathcal{G}(X_1\oplus\cdots\oplus X_n)$;\\[3mm]
$(iii)$ $(a)$ $D_1\in\mathcal{G}_l(X_1)$ and $D_n\in\mathcal{G}_r(X_n)$;\\
\hspace*{8.2mm}$(b)$ $\bigoplus\limits_{s=1}^{i-1}\mathcal{K}(D_s)\prec\mathcal{N}(D_i)$ does not hold for $2\leq i\leq n$ and $\bigoplus\limits_{s=i+1}^n\mathcal{N}(D_s)\prec \mathcal{K}(D_i)$ does not hold for $1\leq i\leq n-1$.\\[3mm]
Then $(i) \Rightarrow (ii).$\\
If $X_1,...,X_n$ are infinite dimensional Hilbert spaces, then $(ii)\Rightarrow(iii)$.
\end{Teorema}
\textbf{Proof. }$(ii)\Rightarrow(iii)$

Let $T_n^d(A)$ be invertible for some $A\in\mathcal{B}_n$. Then $T_n^d(A)$ is both left and right invertible, and so Theorems \ref{LEVIINVERZ'} and \ref{DESNIINVERZ'} yield the desired.

\noindent$(i)\Rightarrow(ii)$

We find $A\in\mathcal{B}_n$ such that $\alpha(T_n^d(A))=0$ and $\mathcal{R}(T_n^d(A))= X_1\oplus\cdots\oplus X_n$.  We choose $A=(A_{ij})_{1\leq i<j\leq n}$ so that $A_{ij}=0$ if $j-i\neq 1$, that is we place all nonzero operators of tuple $A$ on the superdiagonal. It remains to define $A_{ij}$ for $j=i+1$, $1\leq i<n$. First notice that $A_{i,i+1}: X_{i+1}\rightarrow  X_i$. Since all of diagonal entries have closed ranges, we know that $ X_{i+1}=\mathcal{N}(D_{i+1})\oplus \mathcal{M}(D_{i+1})$, $ X_i=\mathcal{K}(D_i)\oplus\mathcal{R}(D_i)$, and we have $\alpha(D_{i+1})=\beta(D_i)$. It follows that there is an invertible $J_{i}:\mathcal{N}(D_{i+1})\rightarrow \mathcal{K}(D_i)$. We put $A_{i,i+1}=\begin{bmatrix}J_{i} & 0\\ 0 & 0\end{bmatrix}:\begin{bmatrix}\mathcal{N}(D_{i+1})\\ \mathcal{M}(D_{i+1})\end{bmatrix}\rightarrow\begin{bmatrix}\mathcal{K}(D_i)\\ \mathcal{R}(D_i)\end{bmatrix}$, and we implement this procedure for all $1\leq i\leq n-1$. 

Notice that $\mathcal{R}(A_{i,i+1})=\mathcal{K}(D_i)$ for each $1\leq i\leq n-1.$ Thus, we prove that $T_n^d(A)$ is surjective in the same way as in the proof of Theorem \ref{DESNIINVERZ'}.

Next, we are able to show that $\mathcal{N}(T_n^d(A))\cong\mathcal{N}(D_1)$, implying $\alpha(T_n^d(A))=\alpha(D_1)=0$. This is proved in the same way as in the proof of Theorem \ref{LEVIINVERZ'}. $\square$
\begin{Posledica}\label{POSLEDICA9'}
Let $D_1\in\mathcal{B}(X_1),\ D_2\in\mathcal{B}(X_2),...,D_n\in\mathcal{B}(X_n)$. Assume that all $D_s-\lambda,\ 2\leq s\leq n-1,\ \lambda\in\mathds{C}$ are regular operators. Then
$$
\bigcap_{A\in\mathcal{B}_n}\sigma(T_n^d(A))\subseteq\\\sigma_{l}(D_1)\cup\sigma_{r}(D_n)\cup\Big(\bigcup_{k=1}^{n-1}\Delta_k'\Big),
$$
where
$$
\Delta_k':=\Big\lbrace\lambda\in\mathds{C}:\ \mathcal{N}(D_{k+1}-\lambda)\cong \mathcal{K}(D_k-\lambda)\ does\ not\ hold\Big\rbrace,\ 1\leq k\leq n-1.
$$
If $X_1,...,X_n$ are infinite dimensional Hilbert spaces, then
$$
\sigma_{l}(D_1)\cup\sigma_{r}(D_n)\cup\Big(\bigcup\limits_{k=2}^{n-1}\Delta_k\Big)\cup\Delta_n\subseteq\\\bigcap_{A\in\mathcal{B}_n}\sigma(T_n^d(A)),
$$
where
$$
\begin{aligned}
\Delta_k=\Big\lbrace\lambda\in\mathds{C}:\ \bigoplus_{s=1}^{k-1}\mathcal{K}(D_s-\lambda)\prec\mathcal{N}(D_k-\lambda)\ holds\Big\rbrace\cup\\
\Big\lbrace\lambda\in\mathds{C}:\ \bigoplus_{s=k+1}^n\mathcal{N}(D_s-\lambda)\prec \mathcal{K}(D_k-\lambda)\ holds\Big\rbrace,\quad 2\leq k\leq n-1,
\end{aligned}
$$
$$
\begin{aligned}
\Delta_n=\Big\lbrace\lambda\in\mathds{C}:\ \bigoplus_{s=1}^{n-1}\mathcal{K}(D_s-\lambda)\prec\mathcal{N}(D_n-\lambda)\ holds\Big\rbrace\cup\\
\Big\lbrace\lambda\in\mathds{C}:\ \bigoplus_{s=2}^n\mathcal{N}(D_s-\lambda)\prec \mathcal{K}(D_1-\lambda)\ holds\Big\rbrace.
\end{aligned}
$$
\end{Posledica}
\begin{Primedba}
Obviously, $\Big(\bigcup\limits_{k=2}^{n-1}\Delta_k\Big)\cup\Delta_n\subseteq\Big(\bigcup\limits_{k=1}^{n-1}\Delta_k'\Big)$ holds.
\end{Primedba}
If we put $n=2$ we get:
\begin{Teorema}\label{INVERZn=2'}
Let $D_1\in\mathcal{B}(X_1),\ D_2\in\mathcal{B}(X_2)$. Consider the following statements:\\
$(i)$   $(a)$ $D_1\in\mathcal{G}_l(X_1)$ and $D_2\in\mathcal{G}_r(X_2)$;\\
\hspace*{5.5mm}$(b)$ $\mathcal{N}(D_{2})\cong \mathcal{K}(D_1)$;\\[3mm]
$(ii)$ There exists $A\in\mathcal{B}_2$ such that $T_2^d(A)\in\mathcal{G}(X_1\oplus X_2)$;\\[3mm]
$(iii)$ $(a)$ $D_1\in\mathcal{G}_l(X_1)$ and $D_2\in\mathcal{G}_r(X_2)$;\\
\hspace*{8.2mm}$(b)$ $\mathcal{K}(D_1)\prec\mathcal{N}(D_2)$ does not hold and $\mathcal{N}(D_2)\prec \mathcal{K}(D_1)$ does not hold.\\[3mm]
Then $(i) \Rightarrow (ii).$\\
If $X_1,X_2$ are infinite dimensional Hilbert spaces, then $(ii)\Rightarrow(iii)$.
\end{Teorema}
\begin{Posledica}\label{POSLEDICA10'}
Let $D_1\in\mathcal{B}(X_1),\ D_2\in\mathcal{B}(X_2)$. Then
$$
\bigcap_{A\in\mathcal{B}_2}\sigma(T_2^d(A))\subseteq\\\sigma_{l}(D_1)\cup\sigma_{r}(D_2)\cup\Delta',
$$
where
$$
\Delta':=\Big\lbrace\lambda\in\mathds{C}:\ \mathcal{N}(D_{2}-\lambda)\cong\mathcal{K}(D_1-\lambda)\ does\ not\ hold\Big\rbrace.
$$
If $X_1,X_2$ are infinite dimensional Hilbert spaces, then
$$
\sigma_{l}(D_1)\cup\sigma_{r}(D_2)\cup\Delta\subseteq\bigcap_{A\in\mathcal{B}_2}\sigma(T_2^d(A)),
$$
where
$$
\begin{aligned}
\Delta:=\Big\lbrace\lambda\in\mathds{C}:\ \mathcal{K}(D_1-\lambda)\prec\mathcal{N}(D_2-\lambda)\ holds\Big\rbrace\cup\\\Big\lbrace\lambda\in\mathds{C}:\ \mathcal{N}(D_2-\lambda)\prec \mathcal{K}(D_1-\lambda)\ holds\Big\rbrace.
\end{aligned}
$$
\end{Posledica}
Theorem \ref{INVERZn=2'} interpreted in the setting of Hilbert spaces is a special case of \cite[Theorem 2]{HAN}. Notice that Han et al. (\cite{HAN}) have proved the equivalence $(i)\Leftrightarrow(ii)$ of Theorem \ref{INVERZn=2'} in arbitrary Banach spaces. Corollary \ref{POSLEDICA10'} recovers a result of Du and Pan (\cite[Theorem 2]{DU}). Notice, however, that in \cite{DU} separability was used, while our statement is separability-free. 

\noindent\textbf{Acknowledgments}\\[3mm]
\hspace*{6mm}This work has been supported by the Ministry of Education, Science and Technological Development of the Republic of Serbia under Grant No.\\ 451-03-68/2022-14/200125\\[3mm]
\noindent\textbf{Data availability statement}\\[3mm]
\hspace*{6mm}All data generated or analysed during this study are included in this published article.


\begin{thebibliography}{W}


\bibitem{CHEN} Chen A., Hai G. {\it Perturbations of the right and left spectra for operator matrices}. J. Operator Theory \textbf{67} (2012), no. 1, 207–214 

\bibitem{OPERATORTHEORY}  Djordjevi\'c D. S., {\it Perturbations of spectra of operator matrices}, J. Oper. Theory \textbf{48}(3), 467-486 (2002) 

\bibitem{KOLUNDZIJA2} Djordjevi\'{c}, D. S., Kolund\v{z}ija M. Z. {\it Generalized invertibility of operator matrices.} Ark. Mat. \textbf{50} (2012), no. 2, 259–267 

\bibitem{KOLUNDZIJA3} Djordjevi\'{c} D. S., Kolund\v{z}ija M. Z. {\it Right and left Fredholm operator matrices.} Bull. Korean Math. Soc. \textbf{50} (2013), no. 3, 1021–1027 

\bibitem{DJORDJEVIC}   Djordjevi\'c D. S., Rako\v cevi\'c V., {\it  Lectures on generalized inverses}, University of Ni\v s,
Faculty of Sciences and Mathematics, Ni\v s, 2008 

\bibitem{DU} Du H. K., Pan J., {\it Perturbation of spectrums of 2×2 operator matrices}, Proc. Amer. Math. Soc. \textbf{121} (1994), no. 3, 761–766 

\bibitem{HAN} Han J. K., Lee H. Y., Lee W. Y., {\it Invertible completions of 2×2 upper triangular operator matrices}, Proc. Amer. Math. Soc. \textbf{128} (2000), no. 1, 119–123 

\bibitem{HUANG} Huang J., Wu X., Chen A., {The point spectrum, residual spectrum and continuous spectrum of upper-triangular operator matrices with given diagonal entries}. Mediterr. J. Math. \textbf{13} (2016), no. 5, 3091–3100.

\bibitem{HWANG} Hwang I. S., Lee W. Y. {\it The boundedness below of 2×2 upper triangular operator matrices}. Integral Equations Operator Theory \textbf{39} (2001), no. 3, 267–276 

\bibitem{STEFAN} Ivkovi\'{c} S., {\it On upper triangular operator $2\times 2$ matrices over $C^*$-algebras.} Filomat \textbf{34} (2020), no. 3, 691–706 

\bibitem{KOLUNDZIJA} Kolund\v{z}ija M., {\it Right invertibility of operator matrices.} Funct. Anal. Approx. Comput. \textbf{2} (2010), no. 1, 1–5 

\bibitem{KINEZI} Li Y., Sun X. H., Du H. K., {\it The intersection of left (right) spectra of 2×2 upper triangular operator matrices}. Linear Algebra Appl. \textbf{418} (2006), no. 1, 112–121 

\bibitem{KINEZI2} Li G., Hai G., Chen A. {\it Generalized Weyl spectrum of upper triangular operator matrices}, Mediterr. J. Math. \textbf{12} (2015), no. 3, 1059–1067.

\bibitem{RAKIC} Raki\'c D., Djordjevi\'c D, {\it A note on topological direct sum of subspaces}, Funct. Anal. Approx. Comput. \textbf{10} (1) (2018), 9-20

\bibitem{SARAJLIJA2} Sarajlija N. {\it Fredholmness and Weylness of block operator matrices}, to appear in Filomat, {available at http://arxiv.org/abs/2108.12425}

\bibitem{WU} Wu X., Huang J., {\it Essential spectrum of upper triangular operator matrices}, Ann. Funct. Anal, \textbf{11} (2020), no. 3, 780–798 

\bibitem{ZGUITTI} Zguitti H. {\it A note on Drazin invertibility for upper triangular block operators}. Mediterr. J. Math. \textbf{10} (2013), no. 3, 1497–1507.

\bibitem{ZHANG} Zhang S., Wu Z., {\it Characterizations of perturbations of spectra of 2×2 upper triangular operator matrices.} J. Math. Anal. Appl. \textbf{392} (2012), no. 2, 103–110 

\bibitem{ZANA} \v{Z}ivkovi\'{c} Zlatanovi\'{c} S. \v{C}., Rako\v{c}evi\'{c} V., Djordjevi\'{c} D. S., {\it Fredholm theory}, University of Ni\v{s}, Faculty of Sciences and Mathematics (to appear)
\end{thebibliography}
\end{document}